\documentclass[12pt]{article}
\usepackage{float}
\usepackage[english]{babel}
\usepackage{amssymb, latexsym, amsmath, amsthm}
\usepackage[left=2cm, right=2cm, top=2cm]{geometry}
\usepackage{float}
\usepackage{tikz-cd}
\usepackage{caption}
\usepackage{float}
\usepackage{subcaption}
\usepackage{enumitem}
\usepackage{verbatim}
\usepackage{blindtext}
\usepackage{titlesec}
\titleformat{\section}
   {\normalfont\fontsize{12pt}{14pt}\selectfont\bfseries}{\thesection}{1em}{\MakeUppercase}
\usepackage{sectsty}
\sectionfont{\centering}
\usepackage{bbm}
\usepackage{longtable}
\usepackage{fullpage} % Package to use full page
\usepackage{parskip} % Package to tweak paragraph skipping
\usepackage{hyperref}
\usepackage{enumitem}
\usepackage{mathtools}
\usepackage{xcolor}
\hypersetup{backref=true,       
    pagebackref=true,               
    hyperindex=true,                
    colorlinks=true,                
    breaklinks=true,                
    urlcolor= black,                
    linkcolor= black,                
    bookmarks=true,                 
    bookmarksopen=false,
    filecolor=black,
    citecolor=black,
    linkbordercolor=red
}

\newcommand{\compconj}[1]{%
  \overline{#1}%
}
\theoremstyle{plain}
\newtheorem{thm}{Theorem}[section]

\theoremstyle{definition}
\newtheorem{defn}[thm]{Definition}

\newtheorem{remark}[thm]{Remark}

\usepackage{url}
\usepackage{geometry}
\usepackage{ragged2e}
\AtBeginDocument{\hypersetup{pdfborder={0 0 1}}}
\begin{document}
\title{Isomorphism of Skew-Holomorphic Harmonic Maass-Jacobi Forms and Certain Weak Harmonic Maass Forms}
\author{Ranveer Kumar Singh}
\date{}
\maketitle
\begin{abstract}
\noindent Recently Bringmann, Raum and Richter generalised the definition of Jacobi forms and Skoruppa's skew-holomorphic Jacobi forms by intertwining with harmonic Maass forms. We prove the isomorphism of the Kohnen's plus space analogue of harmonic Maass forms of weight $k-1/2$ for $\Gamma_0(4m)$ and the  space of these skew-holomorphic harmonic Maass-Jacobi forms of weight $k$ and index $m$ for $k$ odd and $m=1$ or a prime. 
\end{abstract}
\section{Introduction}
Let $k$ be an integer, and $m$ a positive integer. Let $M_{k-\frac{1}{2}}\left(\Gamma_{0}(4 m)\right)$ and $M_{k-\frac{1}{2}}^{!}\left(\Gamma_{0}(4 m)\right)$ be the space of holomorphic and weakly holomorphic modular forms respectively, of weight $k-\frac{1}{2}$ for $\Gamma_{0}(4 m).$ Let $M_{k-\frac{1}{2}}^{+}\left(\Gamma_{0}(4 m)\right)$ be the Kohnen's plus space which is a subspace of $M_{k-\frac{1}{2}}\left(\Gamma_{0}(4 m)\right)$ defined as
\[
M_{k-\frac{1}{2}}^{+}\left(\Gamma_{0}(4 m)\right):=\left\{f \in M_{k-\frac{1}{2}}\left(\Gamma_{0}(4 m)\right) |~c_{f}(n)=0 \text { unless }(-1)^{k} n \equiv 0,1\bmod 4 m\right\}.
\]
Eichler and Zagier systematically developed the theory of Jacobi forms in \cite{DZ}. They also proved that the Kohnen's plus space of holomorphic modular forms of weight $k-1/2$ for $\Gamma_0(4m)$ is isomorphic to the space of holomorphic Jacobi forms of weight $k$ and index $m$ for $k$ even and $m=1$ or a prime. They did so by decomposing the Fourier expansion of Jacobi forms in terms of theta series and showing that the coefficients transform like modular forms. For $k$ odd, Skoruppa showed that the plus space is isomorphic to the space of a new kind of Jacobi form which is real analytic in variable $\tau$ and holomorphic in $z$. These forms were called the skew-holomorphic Jacobi forms \cite{S1,S2}. \par 
In 1920, Ramanujan, in his last letter to Hardy \cite{SR} listed 17 functions which he called mock theta functions. These functions led to huge developments in the theory of automorphic forms. After about eight decades, Zwegers \cite{zw} in 2002 came up with a systematic framework for mock theta functions and linked it to harmonic Maass forms. Motivated by this, the theory of Jacobi forms was generalised to harmonic Maass-Jacobi forms  by Bringmann and Richter \cite{BR}. With this generalisation, Cho defined the analogue of Kohnen's plus space for harmonic Maass forms and proved Zagier type isomorphism between the plus space of harmonic Maass forms and the space of harmonic Maass-Jacobi form \cite{cho}. To be precise they proved the isomorphism of plus space of harmonic Maass forms and the space of vector valued harmonic Maass forms. \par 
Recently the theory of skew holomorphic Jacobi forms was generalised to include harmonic properties with respect to some differential operator by Bringmann, Raum and Richter \cite{BRR}. Several different generalisations and structure results have been proved in \cite{BRR}, but we will only consider a subspace of the generalisation containing Skoruppa's skew-holomorphic Jacobi forms. The aim of the present article is to prove the isomorphism of the plus space of harmonic Maass forms and skew-holomorphic harmonic Maass-Jacobi forms combining the results of \cite{BRR,cho}, in the spirit of Skoruppa. To describe the main results, we first introduce notations and definitions. 
\section{Preliminaries}  
\subsection{Harmonic Maass Forms}
In this section, we will briefly review harmonic Maass forms. The reader is referred to \cite{Ono1} for details.\\
Let $\mathbb{H}$ denote the usual upper half space. We write $\tau\in\mathbb{H}$ as $\tau=u+iv$ and $z\in\mathbb{C}$ as $z=x+iy$. Put $q=e(\tau)=e^{2\pi i\tau}$ and $\zeta=e(z)=e^{2\pi iz}.$ Define   
\begin{equation*}
\Gamma_0(N)\coloneqq \left\{\begin{pmatrix}
a&b\\c&d\\
\end{pmatrix}\in \mathrm{SL}_2(\mathbb{Z}):c~\equiv ~0~(\text{mod} ~N)\right\}.
\end{equation*}
Define the weight-$k$ $(\in\mathbb{R})$ hyperbolic Laplacian 
\begin{equation*}
\Delta_k\coloneqq -v^2\left(\frac{\partial^2}{\partial u^2}+\frac{\partial^2}{\partial v^2}\right)+ikv\left(\frac{\partial}{\partial u}+i\frac{\partial}{\partial v}\right)=-4v^2\frac{\partial}{\partial \tau}\frac{\partial}{\partial \bar\tau}+2ikv\frac{\partial}{\partial \bar\tau}.
\end{equation*}
We define Harmonic Maass forms following Bruiner and Funke \cite{BF}.
\theoremstyle{defn}
\begin{defn}
\normalfont
Let $k\in\frac{1}{2}\mathbb{Z}$. A smooth function (in real sense) $f:\mathbb{H}\rightarrow\mathbb{C}$ is called a weight-$k$ harmonic Maass form on $\Gamma_0(N)$ ($4|N$ if $k\in\frac{1}{2}\mathbb{Z}\setminus\mathbb{Z}$) if 
\begin{enumerate}[label=(\roman*)]
\item For all $\begin{pmatrix}
a&b\\c&d\\
\end{pmatrix}\in \Gamma_0(N)$ and $\tau\in\mathbb{H}$, we have
\[
  f\left(\frac{az+b}{cz+d}\right) =
  \begin{cases}
                                   (cz+d)^kf(\tau) & \text{if $k\in \mathbb{Z}$} \\
                                   (\frac{c}{d})\varepsilon_d^{-2k}(cz+d)^kf(\tau) & \text{if $k\in \frac{1}{2}+\mathbb{Z}.$} \\  
  \end{cases}
\] 
Here $\left(\frac{c}{d}\right)$ is the Jacobi symbol and $\varepsilon_d=\sqrt{\left(\frac{-1}{d}\right)}$. Here $\sqrt{.}$ denotes the principal branch of square root.
\item $\Delta_k(f)=0$.
\item There exists a polynomial $P_f(\tau)\in\mathbb{C}[q^{-1}]$ such that $f(\tau)-P_f(\tau)=O(e^{-\varepsilon v})$ as $v\rightarrow\infty$ for some $\varepsilon>0$. Similar conditions hold at other cusps.
\label{(iii)1.1}
\end{enumerate}
If the third condition in the above definition is replaced by $f(\tau)=O(e^{\varepsilon v})$, then $f$ is said to be a harmonic Maass form of \emph{manageable growth}. Space of harmonic Maass forms of weight $k$ is denoted by $H_k(\Gamma_0(N))$ and that of harmonic Maass forms of manageable growth is denoted by $H_k^{!}(\Gamma_0(N))$. 
\label{def 1.1}
\end{defn}
\begin{remark}
One easily sees that 
$M_{k}^!(\Gamma_0(N))\subset H_{k}(\Gamma_0(N))\subset H_{k}^!(\Gamma_0(N)).$
\end{remark} 
$f\in H_k^{!}(\Gamma_0(N))$ has a Fourier expansion \cite{BF} of the form  
\begin{equation}
f(\tau)=f(u+iv)=\sum\limits_{n>> -\infty}c_f^{+}(n)q^n+ c_f^{-}(0)v^{1-k}+\sum\limits_{\substack{n<< \infty\\n\neq 0}}c_f^{-}(n)\Gamma(1-k,-4\pi nv)q^n.
\label{1}
\end{equation}
where $\Gamma(s,z)$ is the incomplete gamma function defined as 
\begin{equation}
\Gamma(s,z)\coloneqq\int\limits_z^{\infty}e^{-t}t^s\frac{dt}{t}.
\label{2}
\end{equation}
The notation $\sum\limits_{n>> -\infty}$ means $\sum\limits_{n=\alpha_f}^{\infty}$ for some $\alpha_f\in\mathbb{Z}$. $\sum\limits_{n<<\infty}$ is defined similarly. 
If $f \in H_{k}(\Gamma)$ then
\begin{equation}
\qquad f(\tau)=f(u+i v)=\sum_{n>>-\infty} c_{f}^{+}(n) q^{n}+\sum_{n<0} c_{f}^{-}(n) \Gamma(1-k,-4 \pi n v) q^{n}.
\label{3}
\end{equation} 
We call 
\begin{equation*}
f^+(\tau)=\sum\limits_{n>> -\infty}c_f^{+}(n)q^n
\end{equation*}
the \emph{holomorphic part} of $f$ and 
\begin{equation*}
f^-(\tau)=c_f^{-}(0)v^{1-k}+\sum\limits_{\substack{n<< \infty\\n\neq 0}}c_f^{-}(n)\Gamma(1-k,-4\pi nv)q^n
\end{equation*}
 the \emph{nonholomorphic part} of $f$.
 \begin{comment}
Define the Bruiner and Funke pairing map (also called shadow operator) as 
\begin{equation*}
\xi_k\coloneqq 2iv^k\compconj{\frac{\partial}{\partial\bar\tau}}.
\end{equation*} 
It is related to the Laplacian operator by 
\begin{equation}
\Delta_k=-\xi_{2-k}\circ\xi_{k}.
\label{4}
\end{equation}
It is known that for $f\in H_{2-k}^{!}(\Gamma_0(N)),~\xi_{2-k}(f)\in M_{k}^!(\Gamma_0(N))-$the space of weight-$k$ weakly holomorhic modular forms. Furthermore, given the Fourier expansion of $f$ as in Eq. \eqref{1}, we have 
\begin{equation}
\xi_{2-k}(f(\tau))=\xi_{2-k}(f^-(\tau))=(k-1)\compconj{c_f^-(0)}-(4\pi)^{k-1}\sum\limits_{n>>-\infty}\compconj{c_f^-(-n)}n^{k-1}q^n.
\label{5}
\end{equation}
\end{comment}
%The image of $f$ under $\xi$ is called the shadow of $f$. Moreover this map is surjective. A \emph{mock modular form} of weight $2-k$ is the holomorphic part $f^+$ of a harmonic Maass form of weight $2-k$ for which $f^-$ is non trivial. $\xi_{2-k}(f)$ is called the shadow of the mock modular form $f^+$. A general question then is to construct the preimage of a given modular form. An explicit construction of mock modular form whose shadows is $\Theta(\tau)^3$ where $\Theta(\tau)$ is the classical theta function has been done by Rhoades and Waldherr in \cite{RW}. Quite recently Herrero and Pippich constructed mock modular forms whose shadows are Eisenstein series of arbitrary integral weight, level and character \cite{HP}.
\subsection{Vector Valued Harmonic Maass Forms}
Let $\mathrm{Mp}_2(\mathbb{R})$ be the metaplectic two-fold cover of $\mathrm{SL}_2(\mathbb{R})$ consisting of elements of the form $(A,\varphi)$ where $A=\begin{psmallmatrix}
a&b\\c&d
\end{psmallmatrix}\in \mathrm{SL}_2(\mathbb{R})$ and $\varphi:\mathbb{H}\rightarrow\mathbb{C}$ is a holomorphic  function such that $\varphi^2(\tau)=c\tau+d$. The operation in $\mathrm{Mp}_2(\mathbb{R})$ is defined as 
\[
(A,\varphi)(B,\psi)=(AB,\varphi(B\tau)\psi(\tau)).
\]
Let $\mathrm{Mp}_2(\mathbb{Z})$ be the inverse image of $\mathrm{SL}_2(\mathbb{Z})$ under the covering map. One can show that $\mathrm{Mp}_2(\mathbb{Z})$ is generated by $\widetilde{T}:=\left(\begin{psmallmatrix} 1& 1 \\ 0 & 1\end{psmallmatrix}, 1\right)$ and $\widetilde{S}:=\left(\begin{psmallmatrix} 0& -1 \\ 1 & 0\end{psmallmatrix}, \sqrt{\tau}\right)$.

Let $V$ be a rational vector space over $\mathbb{Q}$ with a non-degenerate quadratic form $Q$. Let $(x,y)=Q(x+y)-Q(x)-Q(y)$ be the associated bilinear form of signature $\left(p,q\right).$ Let $L \subset V$
be an even lattice with dual $L^*$. We denote the standard basis elements of the group algebra
$\mathbb{C}\left[L^* / L\right]$ by $\mathfrak{e}_{\gamma}$ for $\gamma \in L^* / L$. Let $\varrho_{L}$ be the Weil representation of $\mathrm{Mp}_{2}(\mathbb{Z})$
$\operatorname{on} \mathbb{C}\left[L^* / L\right],$ defined by
\[
\begin{array}{l}
\varrho_{L}(\widetilde{T})\left(\mathfrak{e}_{\gamma}\right):=e(Q(\gamma)) \mathfrak{e}_{\gamma} \\
\varrho_{L}(\widetilde{S})\left(\mathfrak{e}_{\gamma}\right):=\frac{e\left(\left(q-p\right) / 8\right)}{\sqrt{\left|L^* / L\right|}} \sum\limits_{\delta \in L^* / L} e(-(\gamma, \delta)) \mathfrak{e}_{\delta}.
\end{array}
\]
\begin{defn}
Let $k \in \frac{1}{2} \mathbb{Z}.$ A holomorphic function $f: \mathbb{H} \rightarrow \mathbb{C}\left[L^*/ L\right]$ is called a weakly holomorphic
modular form of weight $k$ and type $\varrho_{L}$ for the group $\mathrm{Mp}_{2}(\mathbb{Z})$ if it satisfies:
\begin{enumerate}
\item $f(M \tau)=\varphi(\tau)^{2 k} \varrho_{L}(M, \varphi) f(\tau)$ for all $(M, \varphi) \in \mathrm{Mp}_{2}(\mathbb{Z})$
\item $f$ has a Fourier expansion of the form
\[
f(\tau)=\sum_{\gamma \in L^* / L} \sum_{n \in \mathbb{Z}+Q(\gamma) \atop n \gg-\infty} c_{f}(\gamma, n) e(n \tau) \mathfrak{e}_{\gamma}.
\]
\end{enumerate}
\end{defn}
The space of these $\mathbb{C}\left[L^* / L\right]$ -valued weakly holomorphic modular forms is denoted by $M_{k, \varrho_{L}}^{!}$. Similarly we define holomorphic modular forms of weight $k$ and type $\varrho_L$ and denote the corresponding space by $M_{k, \varrho_{L}}.$
\begin{defn}
A smooth function $f: \mathbb{H} \rightarrow \mathbb{C}\left[L^* / L\right]$ is called a harmonic Maass form of weight $k$
and type $\varrho_{L}$ for the group $\mathrm{Mp}_{2}(\mathbb{Z})$ if it satisfies:
\begin{enumerate}
\item $f(M \tau)=\varphi(\tau)^{2 k} \varrho_{L}(M, \varphi) f(\tau)$ for all $(M, \varphi) \in \mathrm{M} \mathrm{p}_{2}(\mathbb{Z})$;
\item $\Delta_{k} f=0$;
\item There is a polynomial $P_{f}(\tau)=\sum_{\gamma \in L^* / L} \sum_{n \in \mathbb{Z}+Q(\gamma) \atop-\infty \ll n\leq 0} c_{f}^{+}(\gamma, n) e(n \tau) \mathfrak{e}_{\gamma}$ such that
$f(\tau)=P_{f}(\tau)+O\left(e^{-\varepsilon v}\right)$ as $v \rightarrow \infty$ for some $\varepsilon>0$.
\end{enumerate}
\end{defn}
We denote by $H_{k, \varrho_{L}}$ the space of these $\mathbb{C}\left[L^* / L\right]$-valued harmonic Maass forms. One can similarly define harmonic Maass forms of weight $k$
and type $\varrho_{L}$ of manageable growth by replacing Condition 3 above by $f(\tau)=O(e^{\varepsilon v})$. The space of such forms is denoted by $H_{k,\varrho_L}^{!}.$ 
\begin{remark}
It can easily be seen that $M_{k, \varrho_{L}}^{!} \subset H_{k, \varrho_{L}}\subset H_{k,\varrho_L}^{!}.$
\end{remark}
\begin{comment}
For $f \in H_{k, \varrho_{L}}$, we have a unique decomposition $f=f^{+}+f^{-},$ where
\[
\begin{aligned}
f^{+}(\tau) &=\sum_{\gamma \in L^* / L} \sum_{n \in \mathbb{Z}+Q(\gamma)} c_{f}^{+}(\gamma, n) e(n \tau) \mathfrak{e}_{\gamma} \\
f^{-}(\tau) &=\sum_{\gamma \in L^* / L} \sum_{n \in \mathbb{Z}+Q(\gamma)} c_{f}^{-}(\gamma, n) \Gamma(1-k, 4 \pi|n| y) e(n \tau) \mathfrak{e}_{\gamma}.
\end{aligned}
\]
\end{comment}
Define the following subspace of $H_{k-\frac{1}{2}}^!(\Gamma_0(4m))$
\[
H_{k-\frac{1}{2}}^{!+}(\Gamma_0(4m))\coloneqq \{f\in H_{k-\frac{1}{2}}^!(\Gamma_0(4m))~|~c_f^{\pm}(n)=0~\text{unless}~(-1)^{k}n~\equiv~0,1(\text{mod}~4m)\}.
\]
Similarly, one can define $H^+_{k-\frac{1}{2}}(\Gamma_0(4m))$.
Note that $M_{k-\frac{1}{2}}^+(\Gamma_0(4m))\subset M_{k-\frac{1}{2}}^{!+}(\Gamma_0(4m))\subset H_{k-\frac{1}{2}}^{+}(\Gamma_0(4m))\subset H_{k-\frac{1}{2}}^{!+}(\Gamma_0(4m)).$\par
Now for some $m\in\mathbb{Z}_{> 0}$, suppose $(L^*/L,Q)=(\mathbb{Z}/2m\mathbb{Z},Q)$ where $Q(\gamma)=\gamma^2/4m$ for $\gamma\in L^*/L$. In this case the level of $L$ is $4m$ and we have $p-q\equiv 1(\text{mod}~8).$ One easily sees that for this lattice and quadratic form, the Weil representation is given by 
\[
\begin{array}{l}
\varrho_{m}(\widetilde{T}) \mathfrak{e}_{\ell}:=e\left(\frac{\ell^{2}}{4m}\right) \mathfrak{e}_{\ell} \\
\varrho_{m}(\widetilde{S}) \mathfrak{e}_{\ell}:=\frac{1}{\sqrt{2 i m}} \sum\limits_{\ell^{\prime}(\bmod 2 m)} e\left(-\frac{\ell \ell^{\prime}}{2m}\right) \mathfrak{e}_{\ell^{\prime}}.
\end{array}
\]
With these notations, given an $f \in H_{k-\frac{1}{2}}^{!+}\left(\Gamma_{0}(4 m)\right)$ we define a $\mathbb{C}\left[L^*/ L\right]$-valued function $F=\sum_{\gamma \in \mathbb{Z} / 2 m \mathbb{Z}} F_{\gamma} \mathfrak{e}_{\gamma}$ by
\[
F_{\gamma}(\tau):=\frac{1}{s(\gamma)} \sum_{n \in \mathbb{Z}} c_{f}(n, y / 4 m) q^{n / 4 m}
\]
where \[ c_{f}(n, y):=\begin{cases} c_{f}^{+}(n)+c_{f}^{-}(n) \Gamma\left(\frac{3}{2}-k, 4 \pi|n| y\right)&\text{if $n\neq 0$}\\c_{f}^{+}(n)+c_{f}^{-}(n) v^{\frac{3}{2}-k}&\text{if $n=0$}\\\end{cases}\] and $s(\gamma)=1$ if $\gamma \equiv 0, m \bmod 2 m,$ and 2 otherwise. Cho proved the following Theorem \cite{cho}.
\begin{thm}
If $k$ is odd and $m=1$ or a prime, then the map $f \mapsto F$ defines an isomorphism of $H_{k-\frac{1}{2}}^{!+}\left(\Gamma_{0}(4 m)\right)$ onto $H^!_{k-\frac{1}{2}, \varrho_{L}}$.
\label{thm 1.5}
\end{thm}
\begin{remark}
The isomorphism of Theorem \ref{thm 1.5} restricts to the following isomorphisms (cf. Remark 2,(2) of \cite{cho}):
\begin{enumerate}
\item $H_{k-\frac{1}{2}}^{+}\left(\Gamma_{0}(4 m)\right)\simeq H_{k-\frac{1}{2}, \varrho_{L}}$,
\item $M_{k-\frac{1}{2}}^{!+}\left(\Gamma_{0}(4 m)\right)\simeq M^!_{k-\frac{1}{2}, \varrho_{L}}$,
\item $M_{k-\frac{1}{2}}^{+}\left(\Gamma_{0}(4 m)\right)\simeq M^+_{k-\frac{1}{2}, \varrho_{L}}$.
\end{enumerate}
\label{remark 1.7}
\end{remark}
\subsection{Skew-Holomorphic Harmonic Maass Jacobi Forms}
We mostly follow the notation of \cite{BRR} for this section. Let $k,m\in \mathbb{Z}, m>0$. Let $\Gamma^J=\textrm{SL}_2(\mathbb{Z})\ltimes\mathbb{Z}^2$ be the Jacobi modular group. Define the following skew slash operator: For $A=\left(\begin{psmallmatrix}a&b\\c&d\end{psmallmatrix},(\lambda,\mu)\right)$
\begin{equation}
\begin{array}{l}
\left(\left.\Phi\right|_{k, m} ^{\text{sk}} A\right)(\tau, z):=
\quad \Phi\left(\frac{a \tau+b}{c \tau+d}, \frac{z+\lambda \tau+\mu}{c \tau+d}\right)(c \bar{\tau}+d)^{1-k}|c \tau+d|^{-1} e^{2 \pi i m\left(-\frac{c(z+\lambda \tau+\mu)^{2}}{c \tau+d}+\lambda^{2} \tau+2 \lambda z\right)}
\end{array}
\label{•}
\end{equation}
Define the differential operator - the \emph{skew Casimir element} by
\begin{equation}
\mathcal{C}_{k,m}^{\text{sk}}=-\frac{iv^2}{\pi m}v^{\frac{1}{2}-k}\circ \frac{\partial}{\partial\overline{\tau}}\circ v^{k-\frac{1}{2}}L_{m},
\label{•}
\end{equation}
where 
\[
L_m\coloneqq 8\pi im\frac{\partial}{\partial\tau}-\frac{\partial^2}{\partial z^2} 
\]
is the heat operator.
\theoremstyle{defn}
\begin{defn}
\normalfont
A function $\Phi:\mathbb{H}\times\mathbb{C}\rightarrow\mathbb{C}$ is a skew-holomorphic harmonic Maass-Jacobi form of weight $k$ and
index $m$ if $\Phi$ is real-analytic in $\tau\in\mathbb{H}$ and holomorphic in $z\in \mathbb{C}$, and satisfies the following conditions: 
\begin{enumerate}[label=(\roman*)]
\item For all $A\in \Gamma^J$, $\left(\left.\Phi\right|_{k, m} ^{\text{sk}} A\right)=\Phi$;
\item $\mathcal{C}_{k,m}^{\text{sk}}(\Phi)=0$; 
\item $\varphi(\tau,z)=O(e^{\varepsilon v}e^{2\pi my^2/v})$ as $v\rightarrow\infty$ for some $\varepsilon>0$. 
\label{(iii)}
\end{enumerate}
\end{defn} 
We denote the space of skew-holomorphic harmonic Maass-Jacobi form of weight $k$ and index $m$ by $\widehat{\mathbb{J}}^{\text{sk}}_{k,m}$.
\begin{remark}
Let $J_{k,m}^{!\text{sk}}$ (respectively $J_{k,m}^{\text{sk}}$) denote the space of Skoruppa's skew-holomorphic weak (respectively holomorphic) Jacobi forms of weight $k$ and index $m$. Then only easily sees that $J_{k,m}^{\text{sk}}\subset J_{k,m}^{!\text{sk}}\subset \widehat{\mathbb{J}}^{\text{sk}}_{k,m}.$  
\end{remark} 
Proposition 3.6 of \cite{BRR} along with the growth condition \ref{(iii)} implies that $\Phi\in\widehat{\mathbb{J}}^{\text{sk}}_{k,m}$ has Fourier expansion of the form
\begin{equation}
\begin{split}
\Phi(\tau,z)&=v^{\frac{3}{2}-k} \sum_{n, r \in \mathbb{Z} \atop D=0} c^{0}(n, r) q^{n} \zeta^{r}+\sum_{n, r \in \mathbb{Z} \atop D \ll \infty} c^{+}(n, r) \text{exp}\left(\frac{\pi Dv}{m}\right)q^{n}\zeta^{r}\\&+\sum_{n, r \in \mathbb{Z} \atop D \gg-\infty} c^{-}(n,r)\Gamma\left(\frac{3}{2}-k,\frac{\pi D v}{m}\right) \text{exp}\left(\frac{\pi Dv}{m}\right)q^{n} \zeta^{r}
\end{split}
\label{•}
\end{equation}
where $D=r^2-4mn$. We also consider $\Phi\in\widehat{\mathbb{J}}^{\text{sk}}_{k,m}$ which have Fourier expansion of the form
\begin{equation}
\begin{split}
\Phi(\tau,z)&=\sum_{n, r \in \mathbb{Z} \atop D \ll \infty} c^{+}(n, r) \text{exp}\left(\frac{\pi Dv}{m}\right)q^{n}\zeta^{r}\\&\hspace{3cm}+\sum_{n, r \in \mathbb{Z} \atop D >0} c^{-}(n,r)\Gamma\left(\frac{3}{2}-k,\frac{\pi D v}{m}\right)\text{exp}\left(\frac{\pi Dv}{m}\right)q^{n} \zeta^{r}.
\end{split}
\label{•}
\end{equation}
We denote this subspace by $\widehat{\mathbb{J}}^{\text{sk,cusp}}_{k,m}$.
With the above Definition and notations we have the following Theorem:
\begin{thm}
Let $k$ be odd, and $m=1$ or a prime. Then
\begin{enumerate}
\item $\widehat{\mathbb{J}}_{k, m}^{\text{sk}} \simeq H_{k-\frac{1}{2}}^{!+}\left(\Gamma_{0}(4 m)\right)$,
\label{thm 1.10 (1)}
\item $\widehat{\mathbb{J}}_{k, m}^{\text{sk,cusp}} \simeq H_{k-\frac{1}{2}}^{+}\left(\Gamma_{0}(4 m)\right)$.
\label{thm 1.10 (2)}
\end{enumerate}
\label{thm 1.8}
\end{thm}
\section{Proof of Theorem \ref{thm 1.8}}
\label{sec 4}
We will only prove \ref{thm 1.10 (2)} since \ref{thm 1.10 (1)} is similar. Let $\Phi\in\widehat{\mathbb{J}}_{k, m}^{\text{sk,cusp}}$. Then we have that
\begin{equation*}
\begin{split}
\Phi(\tau,z)&=\sum_{n, r \in \mathbb{Z} \atop D \ll \infty} c^{+}(n, r) \text{exp}\left(\frac{\pi Dv}{m}\right)q^{n}\zeta^{r}\\&\hspace{3cm}+\sum_{n, r \in \mathbb{Z} \atop D >0} c^{-}(n,r)\Gamma\left(\frac{3}{2}-k,\frac{\pi D v}{m}\right)\text{exp}\left(\frac{\pi Dv}{m}\right)q^{n} \zeta^{r}.
\end{split}
\end{equation*}
Using the transformation of $\Phi$ for $A=(I,(\lambda,\mu))$, we have that 
\[
\Phi(\tau, z+\lambda \tau+\mu)=e\left(-m\left(\lambda^{2} \tau+2 \lambda z\right)\right) \Phi(\tau, z).
\]
Using similar arguments as in Theorem 2.2 of \cite{DZ}, one can deduce that if $r^{\prime} \equiv r \bmod 2 m$ and $D^{\prime}=D$ with $D^{\prime}:=r^{\prime 2}-4 n^{\prime} m,$ then
\[
c^{+}\left(n^{\prime}, r^{\prime}\right)=c^{+}(n, r),\quad c^{-}\left(n^{\prime}, r^{\prime}\right)=c^{-}(n, r),\quad \Gamma\left(\frac{3}{2}-k, \frac{\pi D^{\prime} y}{m}\right)=\Gamma\left(\frac{3}{2}-k, \frac{\pi D y}{m}\right)
\]
Hence, we can decompose $\Phi(\tau, z)$ as a linear combination of the theta functions as
\[
\Phi(\tau, z)=\sum_{\ell \in \mathbb{Z} / 2 m \mathbb{Z}} h_{\ell}(\tau) \vartheta_{m, \ell}(\tau, z)
\]
where
\begin{equation*}
\begin{split}
h_{\ell}(\tau)&=\sum_{N \gg-\infty} c^{+}\left(\frac{N+r^{2}}{4 m}, r\right)\text{exp}\left(-\frac{\pi Nv}{m}\right) q^{N / 4 m}\\&+\sum_{N<0} c^{-}\left(\frac{N+r^{2}}{4 m}, r\right)\Gamma\left(\frac{3}{2}-k,-\frac{\pi N y}{m}\right)\text{exp}\left(-\frac{\pi Nv}{m}\right) q^{N / 4 m}
\end{split}
\end{equation*}
with any $r \in \mathbb{Z}, r \equiv \ell \bmod 2 m,$ and
\[
\vartheta_{m, \ell}(\tau, z):=\sum_{r \in Z \atop r \equiv \ell \bmod 2 m} q^{r^{2} / 4 m} \zeta^{r}.
\]
Put $g_{\ell}(\tau)=h_{\ell}(-\overline{\tau})$. Then using similar arguments as in Theorem 5.1 of \cite{DZ}, one can show that the $2m$ tuple $(g_{\ell})_{\ell(\text{mod}~2m)}$ satisfies the transformation property of vector valued harmonic Maass form. It remains to show that $\Delta_{k-\frac{1}{2}}g_{\ell}(\tau)=0$. 
Using the fact that $\vartheta_{m,\ell}$ is annihilated by the heat operator $L_m$, it is easy to check that 
\[
\mathcal{C}^{\text{sk}}_{k,m}(h_{\ell}\vartheta_{m,\ell})=2\left[2iv\left(k-\frac{1}{2}\right)\frac{\partial h_{\ell}}{\partial\tau}+4v^2\frac{\partial^2h_{\ell}}{\partial\tau\overline{\tau}}\right]\vartheta_{m,\ell}.
\]
Thus $\mathcal{C}^{\text{sk}}_{k,m}(\Phi)=0$ gives 
\begin{equation}
2iv\left(k-\frac{1}{2}\right)\frac{\partial h_{\ell}}{\partial\tau}+4v^2\frac{\partial^2h_{\ell}}{\partial\tau\overline{\tau}}=0.
\label{9}
\end{equation} 
for each $\ell$. Next observe that 
\[
\Delta_{k-\frac{1}{2}}g_{\ell}(\tau)=\Delta_{k-\frac{1}{2}}h_{\ell}(-\overline{\tau}).
\]
Using the expression for the hyperbolic Laplacian we have 
\[
\Delta_{k-\frac{1}{2}}h_{\ell}(-\overline{\tau})=-4v^2\frac{\partial^2h_{\ell}(-\overline{\tau})}{\partial \tau\overline{\tau}}+2ikv\frac{\partial h_{\ell}(-\overline{\tau})}{\partial \overline{\tau}}
\]
Substituting $T=-\overline{\tau}$ and using chain rule we get
\[
\Delta_{k-\frac{1}{2}}h_{\ell}(-\overline{\tau})=-4(v^2)\frac{\partial^2 h_{\ell}(T)}{\partial T\overline{T}}+2i\left(k-\frac{1}{2}\right)(-v)\frac{\partial h_{\ell}(T)}{\partial T}.
\]
Observing that $V=$ Im$(-\overline{\tau})=v$, we get 
\[
\Delta_{k-\frac{1}{2}}h_{\ell}(-\overline{\tau})=-\left[4V^2\frac{\partial^2 h_{\ell}(T)}{\partial T\overline{T}}+2i\left(k-\frac{1}{2}\right)V\frac{\partial h_{\ell}(T)}{\partial T}\right]=0
\]
where we used Eq. \eqref{9}. Conversely, given $2m$ tuple $(g_{\ell})_{\ell(\text{mod}~2m)}$ of eigenfunctions of the hyperbolic Laplacian, having growth condition as in \ref{(iii)1.1} of Definition \ref{def 1.1} and satisfying the transformation rule of a vector valued harmonic Maass form, it is easy to show that the function
 \[
\Phi(\tau, z)=\sum_{\ell \in \mathbb{Z} / 2 m \mathbb{Z}} h_{\ell}(\tau) \vartheta_{m, \ell}(\tau, z)
\] 
where $h_{\ell}(\tau)=g_{\ell}(-\bar\tau)$, satisfies the transformation property of skew-holomorphic harmonic Maass-Jacobi form. The annihilation by  the skew Casimir element follows easily by a similar calculation as in Eq. \eqref{9}. Finally the Fourier expansion follows from the Fourier expansion \eqref{3} of $g_{\ell}(\tau)$. It follows that we have an isomorphism
\[
\widehat{\mathbb{J}}_{k, m}^{\text{sk,cusp}}\simeq H_{k-\frac{1}{2},\varrho_L}.
\] 
We conclude by noting that Theorem \ref{thm 1.5} implies that
\[
\widehat{\mathbb{J}}_{k, m}^{\text{sk,cusp}}\simeq H_{k-\frac{1}{2},\varrho_L}\simeq H_{k-\frac{1}{2}}^+(\Gamma_0(4m)).
\] 
\begin{remark}
In view of Remark \ref{remark 1.7}, we see that the isomorphism of Theorem \ref{thm 1.8} restricts to isomorphisms of smaller subspaces. To be precise, we have the following isomorphisms: For $k$ odd and $m=1$ or a prime, we have 
\begin{enumerate}
\item $H_{k-\frac{1}{2}}^{!+}\left(\Gamma_{0}(4 m)\right)\simeq H^!_{k-\frac{1}{2}, \varrho_{L}}\simeq \widehat{\mathbb{J}}_{k, m}^{\text{sk}},$
\item $H_{k-\frac{1}{2}}^{+}\left(\Gamma_{0}(4 m)\right)\simeq H_{k-\frac{1}{2}, \varrho_{L}}\simeq \widehat{\mathbb{J}}_{k, m}^{\text{sk,cusp}},$
\item $M_{k-\frac{1}{2}}^{!+}\left(\Gamma_{0}(4 m)\right)\simeq M^!_{k-\frac{1}{2}, \varrho_{L}}\simeq J_{k,m}^{!\text{sk}}$,
\item $M_{k-\frac{1}{2}}^{+}\left(\Gamma_{0}(4 m)\right)\simeq M_{k-\frac{1}{2}, \varrho_{L}}\simeq J_{k,m}^{\text{sk}}$.
\end{enumerate} 
\end{remark}
\bibliographystyle{plain}
\bibliography{mybib}

\addcontentsline{toc}{chapter}{Bibliography}
\begin{thebibliography}{9}
\bibitem{BF} 
 Jan H. Bruinier and Jens Funke, 
\textit{ On two geometric theta lifts}, \textit{Duke Math. J.}, 125 (2004), no. 1, 45–90.
 
\bibitem{BR} 
K. Bringmann and O. K. Richter,  
\textit{ Zagier-type dualities and lifting maps for harmonic Maass-Jacobi
forms}, Adv. Math., 225 (2010), 2298-2315. 
\bibitem{BRR} 
K. Bringmann, M. Raum and O. K. Richter,  
\textit{Harmonic Maass-Jacobi forms with singularities and a theta-like decomposition}, Trans. AMS, 367(9)
 (2015), 6647 - 6670. 
\begin{comment} 
\bibitem{Ram murty} 
 M. Ram Murty,  
\textit{Problems in Analytic Number Theory},  
Second Edition, Springer.
	
\bibitem{Koblitz}
Neal Koblitz,
\textit{Introduction to Elliptic Curves and Modular Forms},
Graduate Texts in Mathematics, Springer, 1993.
\end{comment}	
\bibitem{Ono1} 
Kathrin Bringmann, Amanda Folsom, Ken Ono, Larry Rolen, 
\textit{Harmonic Maass Forms and Mock Modular Forms: Theory and Applications}, 
American Mathematical Society, 2017.
\bibitem{DZ} 
Martin Eichler and Don Zagier, 
\textit{The Theory of Jacobi Forms}  
 Progress in Mathematics book series, Springer 1985.
\begin{comment} 
\bibitem{DMZ} 
Atish Dabholkar, Sameer Murthy and Don Zagier,
\textit{Quantum Black Holes, Wall Crossing, and Mock Modular Forms} 
arXiv:1208.4074 (2012).
\bibitem{RW}
Rhoades, R. C., Waldherr, M., \textit{A Maass lifting of $\vartheta^3$ and class numbers of real and imaginary quadratic fields.} Math. Res. Lett. 18 (2011), no. 5, 1001–1012.
\bibitem{HP}
S. Herrero, A.V. Pippich, \textit{Mock Modular Forms whose Shadows are Eisenstein Series of Integral Weight}, 	arXiv:1809.05690 [math.NT].
\end{comment}
\bibitem{S1}
N.P. Skoruppa, \textit{Explicit formulas for the Fourier coeffcients of Jacobi and elliptic modular forms,} Invent. Math. 102, no. 3 (1990), 501–520.
\bibitem{S2} 
N.P. Skoruppa, \textit{Developments in the theory of Jacobi forms,} Acad. Sci. USSR, Inst. Appl. Math., Khabarovsk (1990), 167–185.
\begin{comment}
\bibitem{Co}
Henri Cohen, 
\textit{ Sums involving the values at negative integers of $L-$functions of quadratic characters,} Math. Ann. 217, 271-285 (1975).
\bibitem{Sh1}
G. Shimura, \textit{Modular forms of half-integral weight. Modular functions of one variable} I. Lecture
Notes in Math. 320, pp. 57-74. Berlin-Heidelberg-New York: Springer 1973.
\bibitem{Sh2}
G. Shimura, \textit{Modular forms of half integral weight,} Ann. of Math. 97, 440-481 (1973).
\bibitem{DZ1}
Don Zagier, \textit{Nombres de classes et formes modulaires de poids 3/2,} C. R. Acad. Sc. Paris 281
(S6r. A), 883-886 (1975).
\bibitem{miyake}
T. Miyake, \textit{Modular forms}, Springer Monographs in Mathematics, (1976).
\bibitem{DZ1}
Don Zagier, \textit{On the values at negative integers of the zeta-function of a real quadratic field,} L'Ens. Math. (1976).
\end{comment}
\bibitem{cho}
Bumkyo Cho, Youngju Choie, \textit{On Harmonic Weak Maass Forms of Half-Integral Weight}, arXiv:1002.1528v3 [math.NT].
\bibitem{zw}
Sander Zwegers, \textit{Mock Theta Functions,} PhD. Thesis, Universiteit Utrecht, The Netherlands, 2002.
\bibitem{SR}
S. Ramanujan, \textit{The lost notebook and other unpublished papers,} Springer, Berlin, 1988. With an introduction by George E. Andrews.
\end{thebibliography}

\end{document}